\documentclass[12pt]{article}

\usepackage{amsfonts}
\usepackage{amsmath}
\usepackage{amssymb}
\usepackage{amscd}
\usepackage{verbatim}
\usepackage{latexsym}
\font\rurm=wncyr10 scaled \magstep1

\begin{document}

\title{Transcendental $\ell$-adic Galois representations}

\author{Chandrashekhar Khare, University of Utah and 
TIFR, \\
Michael Larsen, Indiana University 
\thanks{partially supported
by NSF Grant DMS-0100537}, \\
\& Ravi Ramakrishna, Cornell University 
\thanks{partially supported
by NSF Grant DMS-0102173 and the AMS 
Centennial Research Fellowship.}}

\date{}

\maketitle
\newcommand{\rhobar}{\overline{\rho}}
\newtheorem{theorem}{Theorem}[section]
\newtheorem{lemma}[theorem]{Lemma}
\newtheorem{prop}[theorem]{Proposition}
\newtheorem{fact}[theorem]{Fact}
\newtheorem{cor}[theorem]{Corollary}
\newtheorem{example}{Example}
\newtheorem{conj}[theorem]{Conjecture}
\newtheorem{definition}[theorem]{Definition}
\newtheorem{quest}[theorem]{Question}
\newtheorem{ack}{Acknowledgemets}
\newenvironment{proof}{\smallskip\noindent{\bf Proof.}\ }{$\Box$\smallskip\goodbreak}

\newcommand{\FL}{ {\mathbb F}_{\ell} }
\newcommand{\FF}{ {\mathbb F} }
\newcommand{\Ad}{{\rm Ad}}
\newcommand{\tr}{{\rm tr}\,}
\newcommand{\comp}{{\scriptstyle\circ}}
\def\sha{{{\textnormal{\rurm{Sh}}}}}
\def\qrhob{\bf Q(\bar{\rho})}
\def\eps{\epsilon}
\def\rhobar{ {\bar {\rho} } }
\def\rhob{ {\bar {\rho} } }
\def\wfq{W({\bf{F_q}})}
\def\ad{Ad^0\bar{\rho}}
\def\adst{ (Ad^0\bar{\rho})^*}
\newcommand{\Galois}{\mathrm{Gal}}
\newcommand{\Gal}{\Galois( \bar{ {\mathbb Q}}/{\mathbb Q})}
\newcommand{\C}{{\mathbb C}}
\newcommand{\Q}{{\mathbb Q}}
\newcommand{\Z}{{\mathbb Z}}
\newcommand{\GL}{\mathrm{GL}}
\newcommand{\SL}{\mathrm{SL}}
\newcommand{\R}{{\mathbb R}}
\newcommand{\F}{{\mathbb F}}

\newcommand{\vep}{{\varepsilon}}

\section{Introduction}
Let $F$ be a number field and $G_F=\Galois(\overline{F}/F)$ be its absolute Galois group.
Let ${\mathbb C}_{\ell}$ be the completion of ${\overline {\Q}_{\ell}}$.
In this paper we study continuous,  {\it transcendental} $\ell$-adic 
Galois representations $\rho:G_F \rightarrow \GL_n({\C}_{\ell})$.
Such representations which arise classically (for example from geometry)
have models over a finite extension of ${\Q}_{\ell}$. By an argument that is part of the ``folklore'',
and that we  have seen  attributed to
Florian Pop, one does not get any ``new representations''
if one considers representations $\rho:G_F \rightarrow \GL_n(\overline{{\Q}}_{\ell})$ 
as these have models
over  some $K$ that is  a finite extension of ${\Q}_{\ell}$. This follows simply from the fact 
that a compact subgroup $C$
of $\GL_n(\overline{\Q}_{\ell})$ in fact lies in $\GL_n(K)$ for 
$K$ a sufficiently large finite extension of ${\Q}_{\ell}$ 
({\it Baire category theorem}: namely by this we know that for some finite
$K'/{\Q}_{\ell}$
the intersection of ${C}$ and $\GL_n(K')$ is an open, and hence finite index, 
subgroup of $C$).  

We describe below (see Theorem
\ref{fractious}) the construction  of  an 
example of a transcendental  
representation, which  
\begin{itemize}\item is semisimple, and therefore unramified at a density one set of places by a result we prove below (see Theorem
\ref{ram}), \item  
at a density one set of unramified places has  characteristic 
polynomials of  Frobenii  defined over $\overline{\Q}_{\ell}$,
 \item  does not have a model over 
$\overline{\Q}_{\ell}$, or equivalently over any finite extension of 
${\Q}_{\ell}$. \end{itemize}
By another result we prove below (see Corollary
\ref{finite}), such representations are necessarily {\it infinitely ramified}. In the course of this paragraph, 
we have also given a description
of the three main results of this work.

The study of $\C_{\ell}$-semilinear representations of
$\Galois(\overline {\Q}_\ell/\Q_\ell)$ has been a central object of study in the subject of $p$-adic (or $\ell$-adic in present notation!)
Hodge theory for more than 30 years. 
The $\C_{\ell}$-linear representations that we study here of global Galois groups have to our knowledge not been studied.
The {\it rationale} that we can offer for the representations considered here,
besides the obvious one that we find them diverting, 
may be summarised as follows.
\begin{enumerate}
\item Transcendental $\ell$-adic representations include the study
of big Galois representations, that have been studied earlier in the work of Hida, Mazur et al, as one has continuous 
embeddings ${\Z}_{\ell}[[X_1,\cdots,X_r]] 
\hookrightarrow {\C}_{\ell}$ (for arbitrary $r$). 
\item As the theory of complex representations of various kinds of groups illustrates, it is natural to study representations 
over complete, algebraically closed fields that are further minimal. Here by minimal we mean
in the sense of ${\C}$ being the complete, algebraic closure of $\R$, the $\infty$-adic completion of its prime field $\Q$, 
or $\C_{\ell}$ being the completion of the algebraic closure of $\Q_{\ell}$, the $\ell$-adic completion  of its prime field
$\Q$.
\item This is a continuation of the point above. Galois representations with values in $\GL_n(\C_{\ell})$ can 
arise as limits of representations that are all defined over $\overline \Q_\ell$.
\end{enumerate}

\noindent{\it Acknowledgements:} We would like to thank Dipendra Prasad and Shankar Sen for
some helpful correspondence. We also thank Brian Conrad for spotting a gap in an earlier proof  of the results
in Section 4.

\section{Cebotarev density and ramification}

In this section, we consider continuous Galois representations taking values in the 
field of fractions $K$
of a valuation ring $V$.  We do not assume that $V$ is discrete 
or complete or that its residue field is finite.
Our goal is to prove that continuous semisimple Galois representations taking 
values in $K$ 
are unramified outside a thin set of primes.  This was shown in the locally compact case in
\cite{Kh-Raj}.  We follow the basic strategy of  \cite{Kh-Raj}.
A crucial part of this paper, the estimate of volumes of ``tubular neighborhoods'' of subvarieties
(following J-P.~Serre \cite{S2}), needs to be replaced by a softer technique.  Our method does not give the Hausdorff dimension as in the paper of Serre, but it does allow us to prove that the measure of
the subset of a compact subgroup $\Gamma$ of $\GL_n(K)$ lying in a sufficiently small tubular neighborhood
of a codimension $\ge 1$ subvariety of the Zariski closure of $\Gamma$ can be made as small as
desired. This also proves  a Cebotarev density theorem for transcendental $\ell$-adic representations (see Theorem \ref{cebotarev}) that
might be of independent interest.

For the convenience of the reader, we begin by recalling that the standard argument for finding an
integral basis for representations of compact groups does not depend on the compactness of $V$.

\begin{lemma}
Every finitely generated torsion-free module over a valuation ring is free.
\end{lemma}

\begin{proof}
Let $V$ be a valuation ring with valuation $v$ and $M$ a finitely generated torsion-free $V$-module.  Let
$m_1,m_2,\ldots,m_k$ be a minimal set of generators of $M$.  We claim they form a basis.
If not, after reindexing there exists a linear combination $a_1 m_1 + \cdots + a_r m_r = 0$, where
$1\le r\le k$, all $a_i$ are non-zero, and $v(a_1)\ge v(a_2)\ge\cdots\ge v(a_r)$.  Then, as $M$
is torsion-free,
$$m_r = -\frac{a_1}{a_r}m_1-\cdots-\frac{a_{r-1}}{a_r}m_{r-1}\in V m_1 +\cdots+V m_{r-1},$$
contrary to the minimality of $\{m_i\}$.
\end{proof}

\begin{lemma}\label{lattice}
Let $V$ be a valuation ring with fraction field $K$, $n$ a positive integer, and $\Gamma$ a compact
subgroup of $\GL_n(K)$.   Then $\Gamma$ can be conjugated within $\GL_n(K)$ 
into a subgroup of $\GL_n(V)$.
\end{lemma}

\begin{proof}
Let $\Gamma^\circ = \Gamma\cap\GL_n(V)$.  Then $\Gamma^\circ$ is an open subgroup of $\Gamma$ and therefore of finite index.  Consider the coset decomposition 
$\Gamma =  \gamma_1\Gamma^\circ \cup\cdots\cup\gamma_m\Gamma^\circ$.  Let
$M = \gamma_1 V^n + \cdots +\gamma_m V^n$.  Then $M$ is a finitely generated $V$-submodule of
$K^n$ (which is torsion-free).  Therefore it is torsion-free, and as it spans $K^n$, isomorphic to $V^n$.
It follows that there exists an element of $\GL_n(K)$ mapping the original basis of $V^n$ to a basis of $M$.
\end{proof}

\subsection{Cebotarev density}

Our goal is to show that subvarieties of 
positive codimension in the Zariski closure of the image of a Galois representation
capture only a density-zero set of Frobenius elements.   We do this by considering ``tubular neighborhoods''
of such subvarieties and showing that as ``radius'' goes to 0, measure goes to 0 as well.  The following
proposition is the key.

\begin{prop}
\label{tubular}
Let $\Gamma$ denote a compact subgroup of $\GL_n(K)$, $\mu$ Haar measure on $\Gamma$, 
$G$ the Zariski closure of $\Gamma$ in
$\GL_n$, and $f$ an element of the coordinate ring of $\GL_n$ over $K$ which does not vanish identically
on any component of $G$.  Then for all $\epsilon > 0$, there exists $\alpha$ in the value group of $V$ such that
$$\mu(\{\gamma\in \Gamma\mid v(f(\gamma)) > \alpha\}) < \epsilon.$$

\end{prop}
\begin{proof}
We say that a $K$-subvariety $X\subset G$ is \emph{thin} if for every
finite subset $\{f_1,\ldots,f_m\}$
of the coordinate ring $A$ of $\GL_n$ such that $V(f_1,\ldots,f_m)\cap G = X$, we have
$$\lim_{\alpha\to\infty}\mu(\{\gamma\in \Gamma\mid \forall i\;v(f_i(\gamma)) > \alpha\}) = 0.$$
%
%
The proposition follows immediately from the more general statement that every subvariety $X$ of $G$ of codimension $\ge 1$ is thin.  We prove this by Noetherian induction; the induction step consists in proving
that $X$ is thin if all of its proper subvarieties are so.  

The base case is the empty variety: $X$ is empty if and only if in the coordinate ring $B = A/I$ of $G$, $(\bar f_1,\ldots,\bar f_m) =  B$ or, in other words, if and only if
there exist $a_1,\ldots,a_m\in A$ such that
$$\bar a_1 \bar f_1+\cdots+\bar a_m\bar f_m = 1.$$
In particular,
$$a_1(\gamma) f_1(\gamma) + \cdots + a_m(\gamma) f_m(\gamma) = 1$$
for all $\gamma\in\Gamma$.  Applying $v$, we obtain
$$\min_{1\le i\le m} v(a_i(\gamma)) + v(f_i(\gamma)) \le 0.$$
On the other hand, the sets 
$$\Gamma_{i,\alpha} := \{\gamma\in\Gamma\mid v(a_i(\gamma))\ge -\alpha\}$$
are open, so by compactness, there exists $\alpha$ in the value group of $V$ with
$v(a_i(\gamma))\ge-\alpha$ for all $i$, and for all $\gamma \in \Gamma$.  Thus, 
$$\{\gamma\in \Gamma\mid \forall i\;v(f_i(\gamma)) > \alpha\} = \emptyset.$$

For the induction step, note that $\Gamma$ acts on $\GL_n$ by left-translation and therefore
acts on $A$.  We write $f^\gamma$ for the image of $f\in A$ by $\gamma\in\Gamma$.
For $x\in X$ and $\gamma_1,\ldots,\gamma_i\in \Gamma$,
the set $\{g\in G\mid gx\in \gamma_1 X\cup\cdots\cup \gamma_i X\}$ has dimension
$\dim X < \dim G$.  As $\Gamma$ is dense in $G$, there exists an infinite sequence
$\gamma_1,\gamma_2,\ldots\in\Gamma$ such that $\gamma_i X\neq \gamma_j X$ for all $i\neq j$.
By the induction hypothesis, $\gamma_j^{-1}\gamma_i X\cap X$ is thin for all $i\neq j$.  
Let $(f_1,\ldots,f_k)$ be an ideal in $A$ such that $V(f_1,\ldots,f_k)\cap G = X$.  Then
$$V(f_1,\ldots,f_k,f_1^{\gamma_j^{-1}\gamma_i},\ldots,f_k^{\gamma_j^{-1}\gamma_i})\cap G 
= \gamma_j^{-1}\gamma_i X\cap X,$$
so for all positive integers $n$ there exists $\alpha$ such that 
$$\mu(\{\gamma\in\Gamma\mid \forall h\le k,\;
v(f_h(\gamma)) > \alpha, v(f_h^{\gamma_j^{-1}\gamma_i}(\gamma)) > \alpha\}) < n^{-2}$$
for all $1\le i < j \le n$.  Translating by $\gamma_j$,
\begin{equation}
\label{intersection-bound}
\mu(\{\gamma\in\Gamma\mid \forall h\le k,\;
v(f_h^{\gamma_j}(\gamma)) > \alpha, v(f_h^{\gamma_i}(\gamma)) > \alpha\}) < n^{-2}.
\end{equation}
Let
$$S_j = \{\gamma\in\Gamma\mid \forall h\le k,\;v(f_h^{\gamma_j}(\gamma)) > \alpha\}.$$
By inclusion-exclusion,
$$1 - \sum_{1\le i\le n} \mu(S_i) + \sum_{1\le i < j \le n}\mu(S_i\cap S_j) \ge 0.$$
As
$$\mu(S_j)= \mu(\{\gamma\in\Gamma\mid \forall h\le k,\;v(f_h(\gamma)) > \alpha\})$$
for all $j$, (\ref{intersection-bound}) implies
$$
\mu(\{\gamma\in\Gamma\mid \forall h\le k,\;v(f_h(\gamma)) > \alpha\}) \le
\frac{1+\binom{n}{2}n^{-2}}{n} < \frac{2}{n},
$$
where $n$ can be taken to be as large as we wish.
The proposition follows by induction.
\end{proof}

As an immediate consequence of Proposition \ref{tubular} above 
and the classical Cebotarev density theorem we have:

\begin{theorem}\label{cebotarev}
Let $F$ be a number field, $\overline F$ an algebraic closure, and $G_F = 
\Galois(\overline F/F)$.
Let $V$ be a valuation ring with fraction field $K$ and 
$\rho\colon G_F\to \GL_n(K)$ be a continuous, finitely ramified representation. Let $X$ be a con\-ju\-ga\-tion-invariant subvariety 
of positive codimension in each component of the Zariski closure of $\rho(G_F)$.
Then the set of primes of $F$ at which the Frobenii under $\rho$ lie in $X$ has Dirichlet density $0$.
\end{theorem}

In the case when 
$\rho$ is finitely ramified, we do not know if  
quantitative refinements of the density 0 result, analogous to
Th\'eor\`eme 10 of \cite{S2}, are true in this general setting. 
In any case our ``soft techniques'' will not yield
such quantitative refinements. The theorem of the next section allows to us to get rid of the finitely ramified hypothesis
in the theorem: but without that assumption quantitative refinements cannot be expected as the last section of [KLR] shows.

\subsection{Ramification}

Here is the main theorem of this section.

\begin{theorem}\label{ram}
Let $F$ be a number field, $\overline F$ an algebraic closure, and 
$G_F = \Galois(\overline  F/F)$.
Let $V$ be a valuation ring with fraction field $K$ where $K$ is complete
of characteristic zero and residue characteristic $\ell$.
Let 
$\rho\colon G_F\to \GL_n(K)$ be a continuous semisimple representation.
Then the set of primes of $F$ at which $\rho$ is ramified has Dirichlet density $0$.
\end{theorem}

To prove the main theorem, we need two simple lemmas.  
\begin{lemma}
\label{e-values}
Let $q>1$ be a positive integer, $L$ an algebraically closed field,  
and $x$ and $y$  elements of $\GL_n(L)$ such that
$xyx^{-1} = y^q$.  Then either $y$ is semisimple and of finite order in $\GL_n(L)$ or
there exist eigenvalues $\lambda_1$ and $\lambda_2$ of $x$ with $\lambda_1/\lambda_2 = q$.
\end{lemma}

\begin{proof}
Let $y = y_s y_u$ be the multiplicative Jordan decomposition.  Then
$$x y_s x^{-1} = y_s^q,\ x y_u x^{-1} = y_u^q.$$
If the characteristic polynomial of $y_s$ is $\prod_{i=1}^n x - \lambda_i$, then
there exists a permutation $\pi\in S_n$ such that
$$\lambda_i^q = \lambda_{\pi(i)},$$
which means that
$$\lambda_i^{\prod_{j=1}^n q^j-1} = 1.$$
Thus either $y$ is semisimple of finite order or $y_u\neq 1$.

In the latter case, $x$ and $y_u$ generate a solvable subgroup and therefore a subgroup of $B(L)$,
where $B$ is the Borel subgroup stabilizing some maximal flag of $L^n$.  Let
$$U_1 = [B,B],\,U_2 = [B,U_1],\,\ldots,\,U_n = \{1\}$$
denote the descending central series.  Thus $y_u\in U_1$ but $y_u\not\in U_n$.  Choose $k$ so that
$y_u\in U_k\setminus U_{k+1}$.  Then $B(L)/U_1(L)$ acts on $U_k/U_{k+1}$, and $\bar y_u$ is an eigenvector
of $\bar x\in B(L)/U_1(L)$ with eigenvalue $q$; as the eigenvalues of the diagonal matrix 
$\sum \lambda_i e_{ii}$ acting on $U_k(L)/U_{k+1}(L)$ are 
$$\lambda_i \lambda_{i+k}^{-1},\ 1\le i\le n-k,$$
the lemma follows.
\end{proof}

\begin{lemma}
\label{roots}
In every valuation ring $V$, $1$ has an open neighborhood in which the only root of unity is $1$ itself.
\end{lemma}

\begin{proof}
For $\alpha > 0$ in the value group of $V$, $v(x-1) \ge\alpha$ implies $v(x^n-1)\ge\alpha$.
It therefore suffices to find $\alpha$ such that for every prime $p$, and every primitive $p$th
root of unity $\zeta_p$, $v(\zeta_p-1) \le \alpha$.  
As $\zeta_p-1$ divides $p$, $v(\zeta_p-1)\le v(p)$.  
There is at most one rational prime $p$ in $V$ for which $v(p) > 0$, so we can take
$$\alpha := \sup_p v(p) < \infty.$$
\end{proof}

\noindent
{\bf Proof of Theorem~\ref{ram}.}
The representation $\rho$ can be wildly ramified at $\wp$ only if $\wp$ has  
the same residue characteristic as $V$; henceforth we will ignore this 
finite set of primes.  If $\wp$ is a prime
of $F$ and $G_\wp$ and $T_\wp$ denote the Galois group of the maximal tamely 
ramified extension of
the completion $F_\wp$ and the tame inertia subgroup respectively, 
we have a short exact sequence
$$0\to T_\wp\to G_\wp\to \hat\Z\to 0.$$
The quotient group $\hat \Z$ is topologically generated by the Frobenius class $\sigma_\wp$, $T_\wp$ is
topologically generated by some non-canonical class $\tau_\wp$, and
\begin{equation}
\label{frob-twist}
\sigma_\wp\tau_\wp\sigma_\wp^{-1} = \tau_\wp^{\Vert\wp\Vert}.
\end{equation}
By Lemma~\ref{e-values} and Lemma~\ref{roots}, 
there exists an open neighborhood $G_{F'}$ of the identity in $G_F$
such that no element of $\rho(G_{F'})$ can have an eigenvalue which is a non-trivial
root of unity.  We exclude henceforth from discussion the finite set of primes $\wp'$ 
in the finite extension $F'$ which are ramified over $F$.  Thus
the natural maps $G_{\wp'}\to G_{\wp}$ restrict to isomorphisms $T_{\wp'}\to T_\wp$.
It follows that $\rho(\tau_\wp)$ is unipotent.  

Let $G$ be the Zariski-closure of $\rho(G_F)$ and $Z$ 
the center of the identity component $G^\circ$.
Choose a faithful $K$-representation 
$G/Z\hookrightarrow \GL_m$.   As $1$ is the only unipotent element in 
$Z(K)$,  $T_\wp$ lies in the kernel of the composition map 
$G_F\to G(K)\to (G/Z)(K)\to \GL_m(K)$ 
if and only if $\rho$ is unramified at $\wp$.  
Replacing $\rho$ by this composition if necessary,
we assume without loss of generality that $G$ has semisimple identity component.

By (\ref{frob-twist}) and Lemma~\ref{e-values}, either $\rho$ is unramified at $\wp$ or
$\rho(\sigma_\wp)$ has two eigenvalues in $\bar K$ whose ratio is $\Vert\wp\Vert$.
Let $\vep$ denote the $\ell$-adic cyclotomic character.
Consider the direct sum 
$$\alpha = \rho\oplus\vep\colon G_F\to G\times\GL_1,$$
and let $H$ denote the Zariski closure of $\alpha(G_F)$ 
in this representation.  Thus
$H\subset G\times\GL_1$ and $H$ projects onto each factor.   
By Goursat's lemma, $H$ is the pullback of the graph of an isomorphism between a quotient of
$G$ and a quotient of $\GL_1$.  Every quotient of $\GL_1$ 
is a torus and $G$ admits no non-trivial
toric quotient, so $H = G\times\GL_1$.
Let $X\subset H$ denote the subvariety of
pairs $(g,c)\in H$ such that $g$ and $gc$ have at least one eigenvalue in common.
For each $g$ there are only finitely many possible 
values of $c$, so $X$ is of codimension $\ge 1$ 
in each component of $H$.  By the construction of $X$, if $\rho$ is ramified at $\wp$, then 
for any choice of
$\sigma_\wp$, $\alpha(\sigma_\wp)$ lies in $X$.  By Proposition~\ref{tubular}
(with $\mu$ Haar measure on $\rho(G_F)$),
for any $\epsilon > 0$, we can find an open and closed 
neighborhood $N_\epsilon$ of $X(K)$ in $H(K)$ 
such that $\mu(\rho(G_F)\cap N_\epsilon ) < \epsilon$.
Let $F''$ be a finite Galois extension of $F$ such that 
$\alpha^{-1}(N_\epsilon)$ is a finite union of
$G_{F''}$ cosets.   Thus the image of 
$\alpha^{-1}(X(K))$ in $\Galois(F''/F)$ has less than $\epsilon |\Galois(F''/F)|$
elements and is a union of conjugacy classes in $\Galois(F''/F)$ 
(since $X$ is a union of conjugacy classes
in $H$).  By the Cebotarev density theorem, there is a set 
of primes $\wp$ of $F$ of Dirichlet density 
at least $1-\epsilon$ such that $\rho$ is unramified at $\wp$.  The theorem follows.
$\square$
\vskip1em
\noindent{\bf Remarks:}   1)
As almost all ramification is unipotent, the proof also 
shows that if $\rho$ has abelian image it is {\it finitely ramified}.
\newline\noindent  2)
Just as in \cite{Kh} we can define the notion of a converging sequence of
$\C_{\ell}$ valued Galois representations $\rho_i:G_F \rightarrow \GL_n(\C_\ell)$
in the residually irreducible case, i.e., $\tr(\rho_i(g))$ is a Cauchy sequence uniformly for $g \in G_F$. 
(It is an immediate consequence of  
\cite{Ca} that for each $n$, the  $\rho_i$'s mod $\ell^n$ are then 
eventually constant, and thus have a limit $\rho$.) Just as in \cite{Kh} one proves by essentially the same method as above
that the Dirichlet density of primes that ramify in any of the representations $\rho_i$
is 0.

\section{An example}
For a finite field ${\mathbb F}_{\ell^d}$, we denote the ring
of Witt vectors of ${\mathbb F}_{\ell^d}$ by
$W({\mathbb F}_{\ell^d})$. We assume $\ell>3$ in this section.

\begin{lemma}\label{semi} Let $\ell >3$ and 
$B_m \subset
\GL_2\left(W({\mathbb F}_{\ell^d})/
l^{m}W({\mathbb F}_{\ell^d})\right)$ be a subgroup whose mod $\ell$
reduction
is
$\GL_2({\mathbb F}_{\ell})$. Then $B_m$
contains (up to conjugation by an element of the form
$I+\ell X \in \GL_2\left(W({\mathbb F}_{\ell^d})/
l^{m}W({\mathbb F}_{\ell^d})\right)$)
an element $\left(\begin{array}{cc} a & 0\\0 & 1 \end{array}\right)$ where
$a \in ({\mathbb Z}/l^m{\mathbb Z})^*$ and $a \not \equiv \pm 1$ mod $\ell$.
\end{lemma}
\begin{proof} 
By hypothesis we know that
the mod $\ell$ reduction of $B_m$ contains an element $h$ of the form $\left(\begin{array}{cc} a & 0\\0 & 1 \end{array}\right)$ where
$a \in ({\mathbb Z}/l^m{\mathbb Z})^*$ and $a \not \equiv \pm 1$ mod $\ell$. Choose any lift $g \in B_m$ of this mod $\ell$ element
$h$.
As the characteristic polynomial of $h$ has distinct roots in ${\mathbb F}_{\ell}$, and as $\ell>2$, 
the element $g$ can be conjugated to a 
diagonal matrix viewed as an element of $\GL_2\left(W({\mathbb F}_{\ell^d})/
l^{m}W({\mathbb F}_{\ell^d})\right)$. By taking $\ell^r$th powers, for large enough $r$, 
we obtain an element that has a 
conjugate in $B_m$ with the desired properties.
\end{proof}
\vskip1em\noindent
{\bf Remark:} We need Lemma~\ref{semi} as an ingredient
for Fact~\ref{disjointness}. In \cite{KLR} our $\rho_m$
always had image $\GL_2({\mathbb Z}/{\ell^m}{\mathbb Z})$ and
it obviously  contained elements
like those in Lemma~\ref{semi}.

\begin{definition}\label{rh}
We say a Galois representation $\rhobar$ satisfies our
{\em running hypotheses} if
\begin{itemize}
\item $\rhob:\Gal \rightarrow \GL_2\left({\mathbb Z}/\ell{\mathbb Z}\right)$
      is surjective with $\ell \geq 5$,
\item $det(\rhobar)=\vep$, the cyclotomic character, and
\item $\rhobar$ is unramified
outside a finite set $S$ of primes that necessarily contains
$\ell$.
\end{itemize}
\end{definition}
In this section we assume,
primarily for simplicity, that {\em all determinants are
the cyclotomic character} $\vep$.
We denote by $\ad$ the Galois module of $2 \times 2$ trace zero matrices
over $\FL$ with  Galois action through $\rhobar$ via conjugation.
We denote by $\adst$ the Cartier dual of $\ad$.
We recall some Definitions, Facts, Lemmas and Propositions from
\cite{KLR}.
The precise \cite{KLR} references are given in each case
parenthetically.

\begin{definition}\label{nicedef} (A slight variant
of Definition $1$ of \cite{KLR})
Suppose $\rhob$ satisfies
our running hypotheses.
We say a prime $q$ is {\em nice} (for $\rhob$)
if
\begin{itemize}
  \item $q$ is {\em not} $\pm 1$ mod $\ell$,
  \item $\rhob$ is unramified at $q$,
  \item  the eigenvalues of $\rhob(\sigma_q)$ (where $\sigma_q$ is Frobenius
         at $q$) have ratio $q$.
\end{itemize}
Let $\rho_m$ be a
deformation of $\rhob$ to
$\GL_2\left(W({\mathbb F}_{\ell^d})/\ell^mW({\mathbb F}_{\ell^d})\right)$
with determinant the cyclotomic character.
We say a  prime $q$ is $\rho_m$-{\em nice} if
\begin{itemize}
  \item $q$ is nice for $\rhob$,
  \item $\rho_m$ is unramified at $q$, and the (necessarily distinct) roots
        of the characteristic polynomial
        of $\rho_m(\sigma_q)$ have ratio $q$.
Note that since $q$ is nice,
the mod $\ell^m$ characteristic polynomial of $\rho_m(\sigma_q)$
has distinct roots that are units; it follows that the eigenvalues
of $\rho_m(\sigma_q)$ are well-defined in $W({\mathbb F}_{\ell^d})/\ell^mW({\mathbb
F}_{\ell^d})$.
\end{itemize}
\end{definition}

\begin{fact}\label{disjointness2}
(Weak version of Fact $5$ of \cite{KLR}) Suppose $\rhob$ satisfies
our running hypotheses.
Let the sets $\{f_1,...,f_n\}$  
and $\{ \phi_1,...,\phi_r \}$
be
linearly independent in
$H^1(\Gal,\ad)$ and
$H^1(\Gal,\adst)$ respectively.
Let $\mathbb{Q}(Ad^0(\rhobar))$ be the field
fixed by the kernel of the action of $G_{\mathbb{Q}}$
on $Ad^0(\rhobar)$. Let ${\bf K}=\mathbb{Q}(Ad^0(\rhobar),\mu_{\ell})$
be the field obtained by adjoining the $\ell$th roots of unity to
$\mathbb{Q}(Ad^0(\rhobar))$.
We denote by ${\bf K}_{f_i}$ and ${\bf K}_{\phi_j}$
the fixed fields of the kernels of the restrictions
of $f_i$ and  $\phi_j$ to $G_{{\bf K}}$,
the absolute Galois group of ${\bf K}$.
Then, as ${\mathbb Z}/\ell{\mathbb Z}[\Galois({\bf K}/{\mathbb Q})]$-modules,
$\Galois({\bf K}_{f_i}/{\bf K})$ and $\Galois({\bf K}_{\phi_j}/{\bf K})$
are isomorphic, respectively, to $\ad$ and $(\ad)^*$ and
each of the fields ${\bf K}_{f_i}$ and ${\bf K}_{\phi_j}$
is linearly disjoint over ${\bf K}$ with the compositum of the others.
Let $I$ be a subset of $\{1,...,n\}$ and $J$ a subset of $\{1,...,r\}$.
Then there exists a Cebotarev set $X$
(a set $X$ of primes of positive density coming
from an application of Cebotarev's Theorem)
of primes $w \not \in S$ such that
\begin{itemize}
\item $w$ is nice,
\item $f_i |_{G_w} \neq 0$  for $i \in I$ and
 $f_i |_{G_w} = 0$ for $i \in \{1,\dots,n\} \backslash I$,
\item $\phi_j|_{G_w} \neq 0$ for $j \in J$ and
$\phi_j|_{G_w} =0$ for $j \in \{1,\dots,r\} \backslash J$.
\end{itemize}
\end{fact}

\begin{fact}\label{Sha2} (Lemma $6$ of \cite{KLR})
Denote by the symbol $\sha^i_X(M)$  the kernel of the localisation
map $H^1(G_X,M) \rightarrow \oplus_{v \in X} H^i(G_v,M)$.
Let $\rhob$ satisfy our running hypotheses. There exists a finite
set $T$ of nice primes such that
$\sha^1_{S \cup T}(\ad)$ and
$\sha^2_{S \cup T}(\ad)$ are trivial. After enlarging $S$, we may
thus assume
$\sha^1_{S}(\ad)$ and
$\sha^2_{S}(\ad)$ are trivial.

\end{fact}
                                                                                    
\noindent{\bf Remark:} 
While Fact~\ref{Sha2} above requires
Fact~\ref{disjointness2}, a study of the proofs
of Lemma~\ref{unram} and Propositions~\ref{cheb} and~\ref{polarisation}
(in \cite{KLR}) shows they
require
only Fact~\ref{disjointness} below.
\vskip1em
\begin{fact}\label{disjointness} (Weak variant of Fact $5$ of \cite{KLR})
Let $\rhob$ satisfy our running hypotheses.
Let $\rho_m$ be a deformation of $\rhob$
to $\GL_2\left(W({\mathbb F}_{\ell^d})/\ell^m W({\mathbb F}_{\ell^d})\right)$
unramified outside $S$ with $det(\rho_m)=\vep$.
Suppose 
$\{\phi_1,...,\phi_r\}$ is independent in $H^1(\Gal,\adst)$.
Let $J$ be a subset of $\{1,...,r\}$.
Then there exists a Cebotarev set $X$
(a set $X$ of primes of positive density coming
from an application of Cebotarev's Theorem)
of primes $w \not \in S$ such that

\begin{itemize}
\item $w$ is $\rho_m$-nice,
\item $\phi_j|_{G_w} \neq 0$ for $j \in J$ and
$\phi_j|_{G_w} =0$ for $j \in \{1,\dots,r\} \backslash J$.
\end{itemize}

\end{fact}
\begin{proof}
By Lemma~\ref{semi} we know 
that the image of $\rho_m$
contains (up to conjugation) an element
an  element $\left(\begin{array}{cc} a & 0 \\ 0& 1\end{array}\right)$ 
with $a \in ({\mathbb Z}/\ell^m{\mathbb Z})^*$ and
$a \not \equiv \pm 1$ mod $\ell$. 
Using the fact that $det(\rho_m)=\vep$, 
one sees from the proof of Fact $5$ of \cite{KLR}
that the Cebotarev set of primes with Frobenius
in the conjugacy class of this element 
and the desired splitting properties in the
fields ${\bf K}_{\phi_i}$ provides the $\rho_m$-nice primes.
\end{proof}
\noindent{\bf Remark:} 
It is not true that one can also include in Fact~\ref{disjointness}
splitting properties with independent elements of 
$H^1(\Gal,\ad)$. Such a statement is true for $\GL_2\left(W( \FF_{\ell^d})
/\ell^m W( \FF_{\ell^d})\right)$ representations provided the minimal field
of definition of the adjoint representation is $\FF_{\ell^d}$. 
In our situation the minimal field of definition of the adjoint
representation is always $\FF_{\ell}$.
\vskip1em

\begin{lemma}\label{unram} (Lemma $8$ of \cite{KLR})
Let $\rhob$ satisfy our running hypotheses.
Let $\rho_m$ be a deformation of $\rhob$
to
$\GL_2\left(W({\mathbb F}_{\ell^d})/\ell^m W({\mathbb F}_{\ell^d})\right)$
unramified outside a set $S$.
Suppose
$\sha^1_S(\ad)$ and $\sha^2_S(\ad)$ are trivial and det$(\rho_m)=\vep$.
Let $R$ be any finite collection
of unramified primes of $\rho_m$ disjoint from $S$.
Then there is a finite set  $Q=\{q_1,...,q_n\}$
of $\rho_m$-nice primes disjoint from $R \cup S$
such that
the maps
\begin{equation}\label{thirteen}
H^1(G_{S \cup R\cup Q},\adst) \rightarrow \oplus_{v \in Q}
H^1(G_v,\adst),\end{equation}
\begin{equation}\label{twelve}
H^1(G_{S \cup R\cup Q},\ad) \rightarrow \oplus_{v \in S \cup R}
H^1(G_v,\ad)\end{equation}
and
\begin{equation}H^1(G_{S \cup Q},\ad) \rightarrow
\left(\oplus_{v \in S}  H^1(G_v,\ad) \right)
\oplus \left(  \oplus_{r\in R} H^1_{nr}(G_r,\ad)
\right)\end{equation} are isomorphisms of ${\mathbb F}_{\ell}$ vector
spaces.
($H^1_{nr}(G_v,M)$ denotes the image of the inflation map
$H^1(G_v/I_v,M^{I_v}) \rightarrow H^1(G_v,M)$).
Upon tensoring with ${\mathbb F}_{\ell^d}$ they can be viewed
as isomorphisms of ${\mathbb F}_{\ell^d}$ vector
spaces.
\end{lemma}
\begin{proof}
That the deformation is to the Witt ring is the
only difference between here and the \cite{KLR} situation.
Fact~\ref{disjointness}  provides us with the necessary $\rho_m$-nice primes for
the proof of the Proposition.
The last statement
about tensoring with ${\mathbb F}_{\ell^d}$ is obvious.
\end{proof}                                                                                    
\begin{prop}\label{cheb} (Proposition $9$ of \cite{KLR})
Following the notation of Lem\-ma~\ref{unram},
let $A$ be any finite set of primes disjoint from
$S \cup R \cup Q$.
Fix $k$ between $1$ and $n$. There exists
a Cebotarev set  $T_k$ of primes $t_k$ such that
\begin{itemize}
\item all $t_k \in T_k$ are $\rho_m$ nice,
\item for any $t_k \in T_k$, the kernel of the map of ${\mathbb F}_{\ell}$
    vector spaces
    \begin{equation} H^1(G_{S \cup Q\cup \{t_k\}},\ad)
        \rightarrow \oplus_{v \in S} H^1(G_v,\ad)
         \oplus_{v \in R} H^1_{nr}(G_v,\ad)\end{equation}
is one dimensional spanned by $f_{t_k}$,
\item $f_{t_k}|_{G_v}=0$ for all $v \in S \cup R
\cup Q \cup A \backslash \{q_k\}$,
\item $f_{t_k}$ is unramified at $G_{q_k}$ and
 $f_{t_k}({\sigma_{q_k}}) \neq 0$.
\end{itemize}
\end{prop}
\begin{proof} 
The only difference between this statement and that of \cite{KLR}
is that the deformation is to the Witt vectors. The proof carries over
word for word.
\end{proof}
\begin{prop}\label{polarisation} (Proposition $10$ of \cite{KLR})
Let the notations be as in Lem\-ma~\ref{unram} and Proposition~\ref{cheb}.
There is a set $\tilde{T}_k$ of one or two primes of $T_k$
such that
\begin{itemize}
  \item there is a linear combination $f_k$ of the elements $f_{t_k}$ for
        $t_k \in \tilde{T}_k$ such that $f_k(\sigma_{t_k})=0$ for
        all $t_k \in \tilde{T}_k$ and $f_k|_{G_{q_k}} \neq 0$,
  \item $j<k$ implies that for
        $t_j \in \tilde{T}_j$ we have $f_{t_k}(\sigma_{t_j})=0$,
  \item $f_k|_{G_v}=0$ for all $v \in S \cup R \cup Q \backslash \{q_k\}$.
\end{itemize}
\end{prop}
\begin{proof}
As in Proposition~\ref{cheb}, the proof 
from \cite{KLR} carries over word for word.
\end{proof}
{\bf Remarks:} 
1) If $\rho_m$ is a deformation of $\rhobar $ to
$\GL_2\left(W({\mathbb F}_{\ell^d})/
l^mW({\mathbb F}_{\ell^d})\right)$ and $f \in H^1(\Gal,\ad \otimes_{\FL}
\FF_{\ell^d})$ then $(I+{\ell}^{m-1}f)\rho_m$ is also a deformation
of $\rhobar $ to
$\GL_2\left(W({\mathbb F}_{\ell^d})/
l^mW({\mathbb F}_{\ell^d})\right)$ that is congruent to $\rho_m$ mod $\ell^{m-1}$.
\newline\noindent
2) In \cite{KLR}  $f_k$ of Proposition~\ref{polarisation}
(Proposition $10$ of \cite{KLR}) spanned
a well defined line (over ${\mathbb F}_{\ell}$).
Given $\rho_m$ to $\GL_2\left({\mathbb Z}/{\ell}^m{\mathbb Z}\right)$,
 we proved in \cite{KLR}
that there was a suitable ${\mathbb F}_{\ell}$
linear combination $h$ of the
of the $f_k$ and an element of $H^1(G_{S \cup Q \cup R},\ad)$
such that $(I+\ell^{m-1}h)\rho_m$
had previously chosen characteristic polynomials at primes of $R$,
and was unobstructed at the primes of $S \cup Q \cup (\cup_k T_k)$.
Here we will use an ${\mathbb F}_{\ell^{2^{m-1}}}$
linear combination.
\newline\noindent
                                                                                    
\begin{theorem}\label{fractious}
There exists a potentially semistable
$2$-dimensional continuous irreducible Galois representation
$\rho\colon G_{\Q} := \Galois(\bar\Q/\Q)\to\GL_2\left(\C_\ell\right)$  ($\ell \geq 5$)
such that the trace of Frobenius
$\tr(\sigma_r)$ belongs to $\bar\Q_\ell\subset \C_\ell$
for all $r$ in a set of primes of $\Q$
of density 1 but $\rho$ is not conjugate to any representation
with values in $\bar\Q_\ell$.
\end{theorem}
                                                                                    
\begin{proof}
The proof follows the general strategy of the proof of Application I of
\cite{KLR}, but in addition needs the refinements above.
                                                                                    
\begin{itemize}
\item Let $\rhob$ satisfy our running hypotheses.
      We call this representation $\rho_1$ from now on.
                                                                                    
\item Use Fact~\ref{Sha2} to  enlarge $S$ to $S_1$
      for which
      $\sha^1_{S_1}(\ad)$ and $\sha^2_{S_1}(\ad)$ are trivial and
      obstructions to lifting can be detected locally.
\item Once and for all choose for each $v \in S_1$ local deformations
      of $\rho_1|_{G_v}$ to $\GL_2\left({\mathbb Z}_{\ell}\right)$. 
      For $v=\ell$ we choose a potentially semistable deformation.
      See \cite{R3} for a proof that these local deformations exist.                                                             
\item Since all local mod $\ell$ representations admit deformations
      to mod $\ell^2$ (see above) and all global obstructions to deforming to
      mod $\ell^2$ can be detected locally, we can
      deform $\rho_1$ to  $\rho_2:\Gal \rightarrow
      \GL_2\left({\mathbb Z}/\ell^2{\mathbb Z}\right)$ with $det(\rho_2)=\vep$.
      We are in a situation where we can apply Lemma~\ref{unram} 
      and Propositions~\ref{cheb} and~\ref{polarisation}.                                                                                    
\item Choose a finite set $R_2$ of unramified primes in $\rho_2$
      beneath some bound, say $b_2$. The actual value of $b_2$ 
      (and our later bounds $b_3,b_4,...$)
      is not important.
      What matters is that $R_2$ contains at least one element
      that we call $r_2$.

      Since all determinants of all deformations will be the cyclotomic character,
      choosing the characteristic polynomial of an unramified prime
      amounts to choosing its trace.
      Thus henceforth we will only speak about choosing the traces of Frobenii.

      Our $\rho_2$ is {\em not} the mod ${\ell}^2$ representation we need,
      so we will alter it by a certain $1$-cohomology class.
      Once and for all we choose 
      traces in $ W({\mathbb F}_{\ell})$
      for all  primes in $R_2 \backslash \{r_2\}$.
      For $r_2$ we choose the trace in $W({\FF}_{\ell^2})\backslash
      W(\FL)$. For all $r \in R_2$  these traces are necessarily   congruent to
      $Trace\left(\rho_2(\sigma_{r})\right)$ mod $\ell$.

\item We  use Lemma~\ref{unram} to adjust $\rho_2$  by a $1$-cohomology
      class $g_2$  with coefficients necessarily in $\ad_{/{\mathbb F}_{\ell}}
      \otimes_{{\FF}_{\ell}} {\FF}_{\ell^2}$ 
      so that the traces of Frobenii of all 
      primes in $R_2$ are (mod $\ell^2$) as previously picked {\em and} so that
      at the places of $S_1$ the local representations are the mod $\ell^2$ reduction
      of the previously chosen characteristic zero local representatons.
      Thus there are no obstructions to deforming to mod $l^3$ at primes of $S_1$.
      Our new representation is $(I+ \ell g_2)\rho_2$.

      The representation $(I+ \ell g_2){\rho}_2$ may be
      ramified at new primes in the finite set $Q_2$ of  Lemma~\ref{unram}.
      There may be obstructions at primes of $Q_2$ to deforming
      to mod ${\ell}^3$.
                                                                                    
\item We use Proposition~\ref{polarisation} to remove obstructions at primes of $Q_2$,
      adjusting by a $1$-cohomology class with coefficients in
      $\ad_{/{\mathbb F}_{\ell}} \otimes_{{\FF}_{\ell}} {\FF}_{\ell^2}$
      while introducing new
      ramified primes in a finite set $T_2$ at which {\em there are no obstructions}
      to deforming to mod $\ell^3$. 
      The trace of Frobenius of all primes in $R_2$ will be as 
      previously chosen. 
      The new map, which we yet again call $\rho_2$, has image
      in $\GL_2\left(W({\mathbb F}_{\ell^2})/\ell^2W({\mathbb F}_{\ell^2})\right)$. 
      Let $S_2=S_1 \cup Q_2 \cup T_2$
      be the set of ramified primes of this new deformation.
      For each $v \in S_2 \backslash S_1$ once and for all choose local
      deformations of (our new) $\rho_2|_{G_v}$ to $\GL_2(W(\FF_{\ell^2}))$.      
      (We have
       no control over whether $\rho_2|_{G_v}$ is ramified. 
       If $\rho_2|_{G_v}$ is unramified, we may choose our characteristic zero
       deformation to be either unramified or ramified as we please.)

      Proposition~\ref{polarisation} (and its proof) guarantees we can do this.
      (See \cite{KLR} for details).
                                           
\item Deform $\rho_2$ to $\rho_3:\Gal \rightarrow
      \GL_2\left(W({\mathbb F}_{\ell^2})/\ell^3W({\mathbb F}_{\ell^2})\right)$.

\item Now fix a set $R_3 \supset R_2$ of unramified primes in
      $\rho_3$ beneath some bound $b_3$. 
      The traces of Frobenii of these primes are
      determined mod $\ell^2$ by $\rho_2$.
      Fix a prime $r_3$ in Let $R_3 \backslash R_2$.

      Our $\rho_3$ is {\em not} the mod ${\ell}^3$ representation we need,
      so we will alter it by a certain $1$-cohomology class.
      Once and for all we choose traces in $ W({\mathbb F}_{\ell^2})$
      for all  primes $r$ in $R_3 \backslash (R_2\cup  \{r_3\})$.
      For $r_3$ we choose the trace in $W({\FF}_{\ell^4})\backslash
      W(\FF_{\ell^2})$. For $r \in R_3 \backslash R_2$
      all these traces are necessarily congruent mod $\ell^2$
      to $Trace(\rho_2\left(\sigma_r)\right)$.

\item  Use Lemma~\ref{unram} to adjust $\rho_3$  by  an $1$-cohomology
       class $g_3$ with coefficients necessarily in 
       $\ad_{/{\mathbb F}_{\ell}} \otimes_{\FL} \FF_{\ell^4}$
       so that the traces of Frobenii of all 
       primes in $R_3$ are (mod $\ell^3$) as  previously
       picked {\em and} so that at the places of 
       $S_2$ the local representations are the mod $\ell^3$ reduction
       of the previously chosen characteristic zero local representatons.
       Thus there are no obstructions to deforming to mod $l^4$ at primes of $S_2$.
       This new representation is $ (I+\ell^2 g_3)\rho_3$.

       The representation $(I+\ell^2 g_3){\rho}_3$ may be
       ramified at new primes in a finte set $Q_3$.
       There may be obstructions at primes of $Q_3$ to deforming
       to mod ${\ell}^4$.

\item Use Proposition~\ref{polarisation}
      to remove obstructions at primes of $Q_3$ while introducing new
      ramified primes in a finite set $T_3$ at which {\em there are no obstructions}
      to deforming to mod $\ell^4$. 
      The trace of Frobenius of all primes in $R_3$ will be as previously
      chosen. Let $S_3=S_2 \cup Q_3 \cup T_3$ be the set of ramified primes
      of this new $\rho_3$. 
      For each $v \in S_3 \backslash S_2$ once and for all choose local
      deformations of (our new) $\rho_3|_{G_v}$ to $\GL_2(W(\FF_{\ell^4}))$.
      (We have
       no control over whether (our new) $\rho_3|_{G_v}$ is ramified.
       If (our new) $\rho_3|_{G_v}$ is unramified, we may choose our characteristic zero
       deformation to be either unramified or ramified as we please.)
       Proposition~\ref{polarisation} (and ist proof) guarantees we can do this.
      (See \cite{KLR} for details).

\item At the $m$th stage we will have a representation
      $$\rho_m:\Gal \rightarrow 
      \GL_2\left( W({\mathbb F}_{\ell^{2^{m-1}}})/
      \ell^mW({\mathbb F}_{\ell^{2^{m-1}}})\right)$$
      that we will be able to deform to 
      $$\rho_{m+1}: \Gal \rightarrow 
      \GL_2\left( W({\mathbb F}_{\ell^{2^{m-1}}})/
      \ell^{m+1}W({\mathbb F}_{\ell^{2^{m-1}}})\right).$$
      We make sure $R_{m+1}\backslash R_m$ has at least one prime $r_{m+1}$ 
      whose trace of Frobenius we choose in $W(\FF_{\ell^{2^m}})\backslash
      W(\FF_{\ell^{2^{m-1}}})$ and
      then using Lemma~\ref{unram} 
      alter $\rho_{m+1}$ in by a suitable 
      $1$-cohomology class $g_{m+1}$
      with coefficients in $\ad \otimes_{\FF} {\mathbb F}_{\ell^{2^{m}}}$.
      Then $(I+ \ell^m g_{m+1})\rho_{m+1}$ 
      will have the trace of Frobenius of all
      primes in $R_{m+1}$  as previously chosen {\em and}
      at all places of $S_m$
      the local representations are the mod $\ell^{m+1}$ reduction
      of the previously chosen characteristic zero local representations.
      Thus there are no obstructions to deforming to mod $l^{m+2}$
      at primes of $S_{m}$.
      This process introduces ramification at new primes and we denote
      the set of these and previously ramified primes by $S_{m+1}$.
      Proposition~\ref{polarisation} allows us,
      once and for all, to choose appropriate characteristic
      zero deformations of $\rho_{m+1}|_{G_v}$ for $v \in S_{m+1} \backslash S_m$.
      We may continue the inductive deformation process.
\end{itemize}

The inverse limit of the compatible system of mod $\ell^m$ representations
is valued in $\GL_2\left(W({\overline {\mathbb F}}_{\ell})\right)$. By the choices
of the traces of the Frobenii of $r_i$ this representation is
{\em not} valued in $\GL_2\left(W(\FF_{\ell^d})\right)$ for any $d$.  
Theorem~\ref{ram} and the Baire category theorem show
this is an example of a representation that
has the properties asserted in the theorem.
\end{proof}

\section{Fields of definition of finitely ramified representations}

We show that in examples like the one above the representation is necessarily
infinitely ramified 
(although in the previous example we could have simply ensured that during the construction).

\begin{theorem}\label{finite'}
Let $G$ be a topologically finitely generated profinite group, 
and $\rho:G \rightarrow \GL_n(\C_\ell)$ a continuous
homomorphism. Let $\{x_i\} \subset G$ be a set of elements such 
that the union of their conjugacy classes is dense in $G$
and such that $\tr(\rho(x_i)) \in \bar \Q_\ell$ for all $x_i \in X$. 
Then there exists a finite extension $L$ of $\Q_\ell$ 
such that $tr(\rho(g))\in L$ for all $g\in G$.
If $\rho$ is semisimple, replacing $L$ by a finite extension if 
necessary, $\rho(G)$ can be conjugated into a subgroup of $\GL_n(L)$.
\end{theorem}

\begin{proof}
By Lemma \ref{lattice} above $\rho$ has a model over 
$\cal O$ the valuation ring of $\C_\ell$, and thus we
consider $\rho$ as taking values in $\GL_n({\cal O})$. By continuity of $\rho$,
for any $m$, $\rho$ mod $\ell^m$ has finite image, and can be regarded
as taking values in $\GL_n({\cal O}_{L_m}/\ell^m{\cal O}_{L_m})$ where 
${\cal O}_{L_m}$ is the ring of integers of a finite extension 
$L_m$ of
$\Q_\ell$.  
This can be all be arranged  because of the surjectivity of the map
$\bar \Z_{\ell} \rightarrow {\cal O}/l^m\cal O$. 
The residual homomorphism that arises from reducing $\rho$
modulo the maximal ideal of $\cal O$
is a homomorphism $\rhobar:G_F \rightarrow \GL_n(k)$ for a $k/\F_\ell$ a finite extension of $\F_\ell$, and we take $k$ to be the minimal such field (i.e.,
it is generated by the finitely many matrix entries of the images of 
$\rho$ mod the maximal ideal of $\cal O$).
We assume for convenience that $L_i$ are Galois over $\Q_\ell$, and also that $L_m \subset L_{m+1}$. 
We also assume that $L_1$ is an unramified extension of $\Q_\ell$ with residue field $k$.

Consider the category ${\cal C}$ of complete local Noetherian $W(k)$-algebras with residue field $k$, and with morphisms in this 
category, local morphisms of $W(k)$-algebras, such that the induced map on residue fields
is the identity. For $A$ in $\cal C$, consider deformations
of $\rhobar$ to $\GL_n(A)$, i.e., continuous homomorphisms
$G \rightarrow \GL_n(A)$ that reduce to  $\rhobar$
modulo the maximal ideal of $A$ 
up to conjugation by matrices that reduce to the identity. 
By the standard theory (see \cite{M}), using that $G$ is finitely generated topologically, there is a versal such
deformation $\rho_{\cal R}:G_F \rightarrow \GL_n(\cal R)$ with $\cal R \in \cal C$. 

We need a lemma:

\begin{lemma}\label{specialisation}
 There is continuous map $\pi:{\cal R} \rightarrow {\cal O}$ of local rings such that
the representation $\pi\comp\rho_{\cal R}$ has the same traces as $\rho$, and in fact
is isomorphic to the chosen integral model of $\rho$.
In the case that $\rhobar$ is residually absolutely irreducible, 
$\pi$ is the unique map with the property that
the representation $\pi\comp\rho_{\cal R}$ has the same traces as $\rho$,
and then automatically $\pi\comp\rho_{\cal R}$ is isomorphic to an integral model of $\rho$ that
in this case is unique. 
\end{lemma}

\begin{proof}
Consider $\rho$ mod $\ell^{s}$ for any $s > 0$,
that is valued in $\GL_n({\cal O}_{L_{s}}/\ell^{s}{\cal O}_{L_{s}})$ 
by what was said above. 
For each positive integer $s$ consider
$A_s=({\cal O}_{L_1}+\ell{\cal O}_{L_2}+\cdots+
\ell^{s-1}{\cal O}_{L_s})/\ell^s{\cal O}_{L_s}$.
Then it is easy to see that $A_{s}$ is an object of $\cal C$ and 
$\rho$ mod $\ell^{s}$ is valued in $\GL_n(A_{s})$, via the natural {\it inclusion} 
$A_s \hookrightarrow {\cal O}/\ell^s{\cal O}$. It is also easy to see by inspection
that if we reduce $A_s$ modulo its ideal 
 $\ell^{s-1}{\cal O}_{L_s}$ we get $A_{s-1}$.
By the versal property of $\cal R$, 
we have a morphism $\pi_{s}:{\cal R} \rightarrow 
A_{s}$
such that $\pi_s\comp\rho_{\cal R}$ is isomorphic to (the chosen integral model of ) $\rho$ mod $\ell^s$.
Here we are using the fact that the image of $\rho$ mod $\ell^s$
is valued in the ring that is the image of the composition of inclusions 
$A_{s} \hookrightarrow {\cal O}_{L_{s}}/\ell^{s}{\cal O}_{L_{s}} \hookrightarrow
{\cal O}/\ell^{s}$ for all $s>0$.
In particular
the traces of $\pi_s\comp\rho_{\cal R}$ and 
$\rho$ mod $\ell^{s}$ coincide. 
At this point, in the case when $\cal R$ is a universal object (as is the case when $\rhobar$
is centralised only by scalars)
we directly see that the 
homomorphisms $\pi_{s}$ form a compatible sequence which gives a homomorphism $\pi:{\cal R} \rightarrow {\cal O}$
such that $\pi\comp\rho_{\cal R}$ is isomorphic to
the chosen integral model of  $\rho$ (and in particular has the same traces as $\rho$).
In the case when $\rhobar$ is residually absolutely irreducible we get the additional
claim of the lemma by the results of \cite{Ca}.

In the case when $\cal R$ is only versal the lemma is slightly more delicate.
We need to observe that as
the rings $A_s$ have finite cardinality and as $\cal R$ is topologically finitely generated, for each $s$ there are only finitely many morphisms
$\alpha_s\colon {\cal R} \rightarrow A_s$---in particular, only finitely many 
with the further property that $\alpha_s\comp\rho_{\cal R}$ is isomorpic to $\rho$ mod $\ell^s$.
By a standard compactness argument, there exists a compatible subsequence 
$\alpha_{t_1},\alpha_{t_2},\ldots$ of homomorphisms, which gives a homomorphism $\pi:{\cal R} \rightarrow {\cal O}$ as before.
\end{proof}

%

\vspace{3mm}

The image ${\cal S}$ of ${\cal R}$ in ${\cal O}$ is a quotient of a complete Noetherian $W(k)$-algebra and is
therefore a complete Noetherian $W(k)$-algebra itself.  We have the following key lemma whose proof we owe to Shankar Sen:

\begin{lemma}
Let $K$ denote the field of fractions of $W(k)$.  Then there exists a finite subextension $L$ of $K\subset\C_\ell$ such that the integral closure of $W(k)$ in ${\cal S}\subset\C_\ell$ is contained in $L$.
\end{lemma}

\begin{proof}
As ${\cal S}\subset\C_\ell$, it is a complete Noetherian integral domain, so
by Cohen structure theory (\cite[0$_{\rm IV}$ 19.8.8]{EGA}), there exists a subring ${\cal S}_0\subset{\cal S}$
such that ${\cal S}_0\cong W(k)[[u_1,\ldots,u_r]]$, and $\cal S$ is a local, module-finite 
${\cal S}_0$-algebra.  Suppose it can be generated by $s$ elements.  If $\cal T$ denotes the integral closure of $W(k)$ in $\cal S$, every element of $t\in \cal T$ is integral over ${\cal S}_0$ and therefore satisfies a minimal monic polynomial equation with coefficients in ${\cal S}_0$ and degree $\le s$.   This polynomial divides the minimal polynomial of $t$ over $W(k)$ and therefore has coefficients in $\bar\Q_\ell\cap W(k)[[u_1,\ldots,u_r]] = W(k)$.  It follows that $t$ lies in an extension of $K$ of degree $\le s$; by Krasner's lemma, the compositum of all such extensions is a finite extension of $K$.
\end{proof}

We go back to the proof of Theorem~\ref{finite'}.  
Let $X$ denote the union of the conjugacy classes of the 
$x_i\in G$.  We know that $\tr(\rho)$ is continuous and maps $X$ into 
$L$, $X$ is dense in $G$, and $L$ is closed in $\C_\ell$.  
It follows that $\tr(\rho)$ maps $G$ into $L$.   If $\rho$ is semisimple,
the theory of pseudorepresentations
(see \cite{T})
shows that after 
replacing $L$ by a suitable finite extension, 
we can conjugate $\rho(G)$ into $\GL_n(L)$.
\end{proof}

\begin{cor}\label{finite}
Let $F$ be a number field with Galois group $G_F = \Galois(\bar F/F)$, 
$\ell$ is a rational prime, and
$\rho\colon G_F\to \GL_n(\C_\ell)$ a continuous, irreducible Galois 
representation ramified at a finite set
of primes.  If there exists a set of primes $\wp$ of $F$ of Dirichlet density 1 such that
the trace of Frobenius $\tr(\rho(\sigma_\wp))$ belongs to $\bar\Q_\ell\subset \C_\ell$, then
$\rho$ is conjugate to a representation with values in a finite extension of $\Q_\ell$.
\end{cor}

\begin{proof} 
This follows from the above theorem, using the following consequence of the Hermite-Minkowski theorem: the Galois group of a pro-$\ell$ extension of a number field $K$, ramified only over
a finite set of primes $S$, is topologically finitely generated.
\end{proof}

\vspace{3mm}

\noindent {\it Addresses:}

\noindent CK: 155 S 1400 E, Dept of Math, Univ of Utah, Salt Lake City, UT 84112, USA
and School of Mathematics, TIFR, Homi Bhabha Road, Mumbai 400 005, INDIA.
e-mail: {\tt shekhar@math.utah.edu, shekhar@math.tifr.res.in}

\vspace{3mm}

\noindent ML:  Department of Mathematics, Indiana University, Bloomington, IN 47405, 
USA. e-mail: {\tt larsen@math.indiana.edu }

\vspace{3mm}

\noindent RR: Department of Mathematics, Cornell University, Malott Hall,
Ithaca, NY 14853, USA. e-mail: {\tt ravi@math.cornell.edu}

\end{document}